\documentclass{paper}
\usepackage{amssymb}
\usepackage[latin1]{inputenc}
\usepackage{graphicx, color}
\usepackage{amsmath,amsthm}
\usepackage[colorlinks]{hyperref}
\usepackage{mathrsfs} 

\newtheorem{theorem}{Theorem}[section]

\newtheorem{corollary}[theorem]{Corollary}
\newtheorem{lemma}[theorem]{Lemma}
\newtheorem{proposition}[theorem]{Proposition}

\theoremstyle{definition}
\newtheorem{definition}[theorem]{Definition}

\theoremstyle{remark}

\newcommand{\RR}{{\mathbb R}}
\newcommand{\NN}{{\mathbb N}}
\newcommand{\ZZ}{{\mathbb Z}}

\newcommand{\lex}{\times_{\operatorname{lex}}}

\newcommand{\el}{$\ell$-}
\newcommand{\Spec}{\operatorname{Spec}}
\newcommand{\Max}{\operatorname{Max}}
\newcommand{\Rad}{\operatorname{Rad}}
\newcommand{\ord}{\operatorname{ord}}
\newcommand{\MV}{\mathcal{MV}}
\newcommand{\MVp}{\mathcal{MV}_P}
\renewcommand{\lg}{\ell\mathcal{G}}
\newcommand{\lgu}{\ell\mathcal{G}_u}

\newcommand{\lmapsto}{\longmapsto}
\renewcommand{\P}{\ensuremath{\mathscr{P}}}

\begin{document}
\title{Representation of Perfect and Local MV-algebras}
\author{
\renewcommand{\thefootnote}{\arabic{footnote}}
\rm Brunella Gerla\footnotemark[1]  \and
\renewcommand{\thefootnote}{\arabic{footnote}}
\rm Ciro Russo\footnotemark[2] \and
\renewcommand{\thefootnote}{\arabic{footnote}}
\rm Luca Spada\footnotemark[2]}
\maketitle
\footnotetext[1]{Dipartimento di Informatica e Comunicazione, Universit\`{a} dell'Insubria. Via Mazzini 5, 21100 Varese, Italy. {\tt brunella.gerla@uninsubria.it} }
\footnotetext[2]{Dipartimento di Matematica ed Informatica. Universit\`{a} di Salerno. Via ponte don Melillo, 84084 Fisciano (SA), Italy}

\begin{abstract}
We describe representation theorems for local and perfect MV-algebras in terms of ultraproducts involving the unit interval $[0,1]$. 
Furthermore, we give a representation of local Abelian \el groups with strong unit as quasi-constant functions on an ultraproduct 
of the reals. All the above theorems are proved to have a uniform version, depending only on the cardinality of the algebra to be embedded, as well as a definable construction in ZFC.\\
The paper contains both known and new results and provides a complete overview of representation theorems for such classes.
\end{abstract}

\section{Introduction}
MV-algebras can be seen, in one of their facets, as a non-idempotent generalization of Boolean algebras.
The lack of idempotency gives to MV-algebras a strong monoidal character.
Nevertheless the typical lattice structure of Boolean algebras can be recovered inside MV-algebras,
by a suitable combination of the primitive connectives.
So, MV-algebras have a double behaviour: on the one hand, they look similar to
monoidal objects such as Abelian groups, on the other, they maintain a strong lattice structure (distributive and bounded), which makes many techniques of lattice theory be appropriate in their study.  Such a twin complexion of MV-algebras is endorsed by Mundici's celebrated natural equivalence between MV-algebras and lattice ordered Abelian groups with strong unit.

Boolean algebras work as the \emph{equivalent algebraic semantics} of classical logic.
Boolean algebras also have a neat and intuitive depiction: modulo isomorphisms, any of
them is an algebra of \emph{special} subsets of some set w.r.t. the elementary operations of union, intersection and complement.
This triple scenario is nowadays well understood, the logical conjunction, the set intersection and the Boolean meet are just different personifications of the same concept. The \emph{special} property mentioned above is formalized through the topology of Stone spaces and allows to \emph{select} the right objects in the full powerset of some set.

MV-algebras are the equivalent algebraic semantics of {\L}ukasiewicz logic, one of the longest-known many-valued logic. Nevertheless the monoidal character of such systems
still escapes a complete understanding.  One of the major achievements in this direction, was made by Di Nola who proved that any MV-algebra can be seen, up to isomorphisms, as a subalgebra of the algebra of functions $^{*}[0,1]^{X}$, where $^{*}[0,1]$ is some ultrapower of the real interval $[0,1]$ and $X$ is some suitable set. The structure carried by $^{*}[0,1]^{X}$ is defined pointwise upon the \emph{standard} MV-algebraic structure of $[0,1]$.
Di Nola's result is in complete accordance with the Stone representation Theorem for
Boolean algebras cited above, for subsets of a set $X$ can be seen as discrete functions in $\{0,1\}^{X}$.
Two specialisations of the above result are worth to be mentioned. In \cite{belsemisimple} Belluce proved that every semisimple MV-algebra
can be embedded in an algebra of functions from a set $X$ in $[0,1]$.
MV-algebras that are semisimple and linearly ordered admit a representation as subalgebras of $[0,1]$.

Under this perspective the above representation for MV-algebras suggests two things: it makes sense to consider MV-algebras as appropriate algebraic structures for \emph{fuzzy sets}
; the \emph{infinitesimal} elements arising from the ultrapower construction are some sort of novelty which is produced by the monoidal structure inside MV-algebras and cannot be easily removed.

Concrete representations, as the aforementioned ones, are extremely useful in the study of abstract algebraic structures. Yet, to take full advantage of such correspondences, a full-flagged duality is called for.  Unfortunately, Di Nola representation remains so far an ``embedding theorem'', or, in categorical terms, a faithful and dense functor.
Finding the set of conditions which enable to invert such a functor is probably one of the most interesting open problems regarding MV-algebras.

In this paper we revise some known results on \emph{non-standard} representation for MV-algebras and some subclasses and show some new ones, together with their corresponding results holding for Abelian $\ell u$-groups. In particular we study the subclasses of \emph{local} and \emph{perfect} MV-algebras and prove a uniform non-standard representation theorem for both of them. Our aim is to present a wide perspective of representations for this classes.

More in details, we show that any perfect MV-algebra can be embedded in an algebra of functions from the \emph{spectrum} of the algebra to an ultrapower of $\Gamma(\ZZ \lex \RR, (1,0))$ (\autoref{perflex}).
After recalling a known result on the representation of  local MV-algebras, we use it to prove a representation theorem for local Abelian $\ell u$-groups as the groups of  \emph{quasi-constant} functions over an ultrapower of $\mathbb{R}$ (\autoref{groupquasiconstant}).
For both theorems we also give a uniform version parametrised only by an upper-bound on the cardinality of the algebras to be embedded (\autoref{uniformlocal}).

\section{Basics on MV-algebras}

\begin{definition}
An MV-algebra is an algebra $\langle A,\oplus ,\lnot ,0\rangle $ with a binary operation $\oplus$, a unary operation $\lnot $ and a constant $0$ satisfying the following equations:

\begin{enumerate}
\item $x\oplus \left( y\oplus z\right) =(x\oplus y)\oplus z$
\item $x\oplus y=y\oplus x$
\item $x\oplus 0=x$
\item $\lnot \lnot x=x$
\item  $x\oplus \lnot 0=\lnot 0$
\item $\lnot \left( \lnot x\oplus y\right)
\oplus y=\lnot \left( \lnot y\oplus x\right) \oplus x.$
\end{enumerate}
\end{definition}


The real unit interval $[0,1]$ equipped with operations $x\oplus
y=\min (1,x+y)$ and $\lnot x=1-x$ is an MV-algebra, denoted henceforth by $[0,1]$. The MV-algebra $[0,1]$ generates the whole variety of MV-algebras and hence it has a prominent role in this theory.

%
Given an MV-algebra $A$ and a set $X$, the set $A^{X}$ of all functions $%
f:X\rightarrow A$ becomes an MV-algebra if the operations $\oplus $ and $%
\lnot $ and the element $0$ are defined pointwise.
In particular the set of all functions from a set $X$ to $[0,1]$
can be equipped with a structure of MV-algebra. We shall see that
MV-algebras of this kind can be characterized in terms of their ideals.

\medskip
A partially ordered Abelian group \cite{cdm} is an Abelian group
$\langle G,+,-,0\rangle $ endowed with a partial order relation
$\leq $ that is compatible with addition. When the order of $G$
defines a lattice structure, $G$ is called a \emph{lattice ordered
Abelian group} or $\ell$-group. A {\em strong unit} of an
$\ell$-group $G$ is an element $u \in G$ such that $u \geq 0$ and
for each $x \in G$ there is an integer $n \geq 0$ with $|x| \leq
nu$.
In the following by \el groups we shall always mean Abelian \el groups
and by $\ell u$-groups we shall mean Abelian \el groups with strong unit.

Let $G$ be an $\ell$-group. For any element $u\in G$, $u>0$, (not
necessarily $u $ being a strong unit of $G$) we let $[0,u] = \{x \in G : 0 \leq x \leq u\}$, and for each $x,y \in [0,u]$, $x\oplus y=u\wedge ( x+y)$ and $\neg x=u-x$.
The structure $\langle [0,u],\oplus, \lnot, 0\rangle$ denoted by $\Gamma (G,u)$, is an MV-algebra.

\medskip

Let $\MV$ denote the category of MV-algebras. 
\begin{theorem}\label{thcdm}{\rm \cite{cdm}}
Let $\lgu$ denote the category whose objects are pairs
$(G,u)$ with $G$ an $\ell u$-group and $u$ a
distinguished strong unit of $G$, and whose morphisms are unital
$\ell$-homomorphisms. Then $\Gamma $
is a natural equivalence (i.e., a full, faithful, dense functor) between
$\lgu$ and  $\MV$.
\end{theorem}
The inverse equivalence $\MV \longrightarrow \lgu$ is usually denoted by $\Xi$.

\autoref{thcdm} allows to describe more examples of MV-algebras. Let $\ZZ$ be the $\ell$-group of integer numbers with natural order and denote by $\ZZ \lex \ZZ$ the $\ell$-group given by the lexicographic product (the order is hence total). The MV-algebra $\Gamma(\ZZ \lex \ZZ, (n-1,0))$ is called {\em Komori chain of rank $n$} and is denoted by $S_n^\omega$ (\cite{cdm},\cite{komori}).
The MV-algebra $S_2^\omega$ is also known as {\em Chang MV-algebra} and it is denoted by $C$ (\cite{chang}).

On each MV-algebra $A$ we define
\[1:=\lnot 0,\quad x\odot y:=\lnot \left( \lnot x\oplus \lnot y\right),\quad x\ominus y:=x\odot \lnot y.\]
Note that in the MV-algebra $[0,1]$ we have $x\odot y=\max \left( 0,x+y-1\right) $ and $x\ominus y=\max (0,x-y).$

Let $A$ be an MV-algebra. For any two elements $x$ and $y$ of $A$
we write $x\leq y $ iff $x\ominus y=0$. It follows that $\leq $
is a lattice order in which the operations are given by
\begin{eqnarray}
\label{eq:lattice}
x\vee y&=&\left( x\odot \lnot y\right) \oplus y=\left(
x\ominus y\right) \oplus y, \\
\label{eq:lattice2}x\wedge y&=&\lnot \left( \lnot x\vee \lnot y\right) =x\odot \left(
\lnot x\oplus y\right) .
\end{eqnarray}
In the algebra $[0,1]$, such an order coincides with the natural one.  An MV-algebra whose  order is total is called an MV-chain.

\begin{definition}
An ideal of an MV-algebra $A$ is a subset $I$ of $A$ such that
$0\in I$, if $x, y\in I$ then $x\oplus y\in I$ and
 if $x\in I$, $y\in A$ and $y\leq x$ then $y\in I$.
\end{definition}

  Any ideal $I$ of $A$ induces a congruence $\equiv _{I}$
on $A$ such that for each $x$ and $y\in A$ we have:
\begin{equation*}
x\equiv _{I}y\Longleftrightarrow d(x,y)=\left( x\ominus y\right)
\oplus \left( y\ominus x\right) \in I.
\end{equation*}
  We indicate with $x/I$ or $[x]_I$  the equivalence class
of $x$ and with $A/I$ the set of all equivalence classes with
respect to $\equiv _{I}.$ The quotient MV-algebra is defined as
usual.

  An ideal of an MV-algebra $A$ is \emph{proper} if
$I\neq A$. We say that $I$ is \emph{prime} iff it is proper and
for each $x$ and $y$ in $A$, either $\left( x\ominus y\right) \in
I$ or $\left( y\ominus x\right) \in I$, hence $A/I$ is an MV-chain.

  An ideal of an MV-algebra $A$ is called \emph{maximal}
if it is proper and it is not strictly contained in any proper ideal.

  Let us denote by $\Spec A$ the set of prime ideals and $\Max A$ the set of maximal ideals in $A$. $\Spec A$ is called the
prime spectrum of $A$ and $\Max A$ the maximal spectrum of $A$.
%
%
%

\begin{theorem}[Chang Subdirect Representation Theorem, \cite{chang59}]\label{theo:chan}
Every nontrivial MV-algebra is a subdirect product of MV-chains. In particular, any MV-algebra $A$ is embeddable in $\prod_{P \in \Spec A}A/P$.
\end{theorem}

By \autoref{theo:chan}, any MV-algebra $A$ is a subdirect product of $\{A/P \mid P \in \Spec A\}$. By
finding a common embedding for each $A/P$,  Di Nola showed the
following representation theorem.

\begin{theorem}{\rm \cite{din1}}\label{theo:DNLT}
Up to isomorphisms, every MV-algebra $A$ is an algebra of $[ 0,1]^\ast$-valued functions over $\Spec A$, where $[0,1]^{\ast }$ is an ultrapower of $[0,1]$, only depending on the cardinality of $A$. 
\end{theorem}

Recently, such a representation have been sharpened to a uniform version that asserts the existence, for any infinite cardinal $\alpha$, of a single MV-algebra of functions in $[0,1]^{\ast}$, which contains all MV-algebras of cardinality less than or equal to $\alpha$.
\begin{theorem}{\rm \cite{DNLS}}
For any infinite cardinal $\alpha$ there exists an MV-algebra of functions $A$ such that any MV-algebra of infinite cardinality at most $\alpha$ embeds in $A$.
\end{theorem}

Such a result is obtained by using standard techniques in model theory, involving $\alpha$-regular ultrafilter \cite[Section 4.3]{ChK}. A further application of known constructions in model theory (iterated ultrapowers \cite[Section 6.5]{ChK}) yields a sort of \emph{canonical} MV-algebra of values, parametrised only by some cardinal number $\alpha$.

\begin{theorem}{\rm \cite{DNLS}}
For every infinite cardinal $\alpha$ there is an iterated ultrapower $\Pi_\alpha$ of $[0,1]$, definable\footnote{By \emph{definable} we mean definable in ZFC} in $\alpha$, such that every MV-algebra of cardinality $\alpha$ embeds in an algebra of functions with values in $\Pi_\alpha$.
\end{theorem}

\section{Perfect MV-algebras}

An {\em infinitesimal} element in an MV-algebra is an element $x\neq 0$ such that $nx < \neg x$,
for every $n \in \NN$. Clearly in the MV-algebra $[0,1]$ and in all MV-algebras of functions
taking values in $[0,1]$, there are no such elements. On the other hand,
there are  MV-algebras that are generated by their infinitesimal elements, as for example
the Chang MV-algebra $S^2_\omega$. MV-algebras generated by their infinitesimals
are called {\em perfect} and they form a subcategory of $\mathcal{MV}$
that is equivalent to the category of \el groups.
The class of perfect MV-algebras is included in the variety generated by the Chang MV-algebra,
but it is not a variety: it is a universal class \cite{beldingerper}.
Infinitesimal elements and perfect MV-algebras are strongly related with the phenomenon of incompleteness
of first order {\L}ukasiewicz logic: indeed the subalgebra of the Lindenbaum algebra of first order {\L}ukasiewicz logic
generated by the classes of formulas which are valid but non-provable is a perfect MV-algebra (see \cite{beldin}, \cite{scarp}).

\begin{definition}
For any MV-algebra $A$, the \emph{radical} of $A$ (denoted by $\Rad A$)
is the intersection of all maximal ideals of $A$.
\end{definition}
All the non-zero elements of the radical of an MV-algebra are infinitesimal.


\begin{definition}
\label{def:simplesemisimple}
An MV-algebra is called {\em simple} if it has exactly two ideals, in other words if it is non trivial and $\{0\}$ is the only proper ideal.
An MV-algebra $A$ is said to be {\em semisimple} if $A$ is non trivial and $\Rad A = \{0\} .$
\end{definition}

Every simple MV-algebra is semisimple. The proof of the following proposition can be
found in \cite{belsemisimple,chang,chang59,cdm}.
\begin{proposition}
An MV-algebra is simple if and only if it is isomorphic to a subalgebra of $[0,1]$. An MV-algebra is semisimple if and only if it is isomorphic to a separating MV-algebra of $[0,1]$-valued continuous functions on some nonempty compact Hausdorff space, with pointwise operations.
\end{proposition}

\begin{definition}
An MV-algebra $A$ is called {\em perfect} if it is nontrivial and
\begin{equation*}
A=\Rad A\cup \lnot \Rad A,
\end{equation*}
  where $\lnot \Rad A=\left\{ x\in A:\lnot x\in
\Rad A\right\} .$
\end{definition}


Let $\MVp$ be the full subcategory of $\MV$ having as objects the perfect MV-algebras.
Following \cite{dinletPerf}, we can associate each perfect MV-algebra $A$ with an $\ell$-group
$G = \mathcal D(A)$ as follows. Let $\theta \subseteq \Rad A^2 \times \Rad A^2$ be defined by $(x,y) \theta (x',y')$ if and only if $x \oplus y' = x'\oplus y$; the relation $\theta$ is a congruence and
$$\mathcal{D}(A) = \langle \Rad A^2/\theta, +, \leq, [0,0]\rangle$$
is an $\ell$-group, with $-[x,y]=[y,x]$. It can be shown that this construction yields a functor between $\MVp$ and the category $\lg$ of $\ell$-groups. Moreover, the functor $\mathcal D$ is a categorical equivalence whose inverse is given by
$$\mathcal G: \ G \in \lg \ \longmapsto \ \Gamma(\ZZ \lex G, (1,0)) \in \MVp.$$

In particular, we have that for every perfect MV-algebra $A$ there exists an \el group $G$ such that
$A$ is isomorphic to $\Gamma(\ZZ \lex G,(1,0))$.

\begin{theorem}{\rm \cite{beldinger}}\label{th:MVultrapower}
Let $\prod_F(G_i)$ denote the ultraproduct of the $\ell$-groups $(G_i)_{i \in I}$ with respect to a non-principal ultrafilter $F$ of $2^I$ and $\prod_F(\mathcal G(G_i))$ denote the ultraproduct of the perfect MV-algebras $\mathcal G(G_i)$ with respect to $F$. Then $\mathcal G(\prod_F(G_i)) \cong\, \prod_F(\mathcal G(G_i))$.
\end{theorem}

\begin{theorem}\label{perflex}
Every perfect MV-algebra can be embedded into an algebra of functions taking values in
an ultrapower of $\Gamma(\ZZ \lex \RR, (1,0))$.
\end{theorem}
\begin{proof}
If $A$ is a perfect MV-chain, then $\mathcal D(A)$ is a totally ordered group,
hence $\mathcal D(A)$ can be embedded in some ultrapower ${}^*\RR$ of $\RR$. Since $\mathcal G$ preserves embeddings, then $\mathcal G(\mathcal D(A))$ embeds in $\mathcal G({}^*\RR)$ and, since $\mathcal{G}$ and $\mathcal D$ form an equivalence between the categories of perfect MV-algebras and $\ell$-groups, we have an embedding of $A$ into $\mathcal{G}({}^*\RR) = \Gamma(\ZZ \lex {}^*\RR,(1,0))$.

From \autoref{th:MVultrapower} it follows that every perfect MV-chain can be embedded into an ultrapower of the perfect MV-chain $\Gamma(\ZZ \lex \RR, (1,0))$ (see also \cite[Theorem 23]{beldinger}).

Let $\mathbf{K}$ be the class of all ultrapowers of $\Gamma(\ZZ \lex \RR,(1,0))$ and $M$ be a perfect MV-algebra. By Chang representation theorem, $M$ is embeddable in the product
$\prod_{P \in \Spec M}M/P$ where each $M/P$ is a perfect MV-chain. Then each $M/P$ can be embedded into an element $N_P$ of $\mathbf{K}$. On the other hand, all algebras in $\mathbf{K}$ are elementarily equivalent hence, by the joint embedding property, there exists a perfect MV-algebra $N \in \mathbf{K}$ such that each $N_P$ embeds in $N$. Hence $A$ embeds in $\prod_{P \in \Spec M}N \cong N^{\Spec M}$ and $N$ is an ultrapower of $\Gamma(\ZZ \lex \RR,(1,0))$.
\end{proof}

We aim now at giving a uniform version of the previous theorem.  The idea is basically the same of \cite{DNLS}, we only have to check that the construction goes through also in the restricted case of perfect MV-algebras.

\begin{definition}
Let $\alpha$ be a cardinal.  A proper filter $D$ over $I$ is said to be \emph{$\alpha$-regular} if there exists a set $E\subseteq D$ such that $|E|=\alpha$ and each $i\in I$ belongs to only finitely many $e\in E$.
\end{definition}

For any set $I$ of infinite cardinality $\alpha$ there exists an $\alpha$-regular ultrafilter over $I$ \cite{ChK}. Given a cardinal $\alpha$, let  $\alpha^{+}$ be the smallest cardinal greater than $\alpha$. We let $\equiv$ stand for elementary equivalence, $\hookrightarrow$ for embeddability and $\hookrightarrow_{\mathit{el}}$ for elementary embeddability.

\begin{definition}
Given a cardinal $\alpha$, we say that a model $\mathfrak{A}$ is \emph{$\alpha$-universal} if and only if for every model $\mathfrak{B}$ we have:
\[ \mathfrak{B}\equiv \mathfrak{A}\quad\text{ and }\quad |\mathfrak{B}| <\alpha\quad \text{ implies }\quad\mathfrak{B}\hookrightarrow_{\mathit{el}}\mathfrak{A}.\]
\end{definition}
\begin{theorem}{\rm \cite{ChK}}\label{thm:chang}
Let $\mathcal{L}$ be a first order language such that $|\mathcal{L}|\leq \alpha$ and $D$ be an ultrafilter which is $\alpha$-regular. Then, for every model $\mathfrak{A}$, the ultrapower $\prod_D \mathfrak{A}$ is $\alpha^+$-universal.
\end{theorem}

Let $[0,1]_{\rm lex}$ be the perfect MV-algebra $\Gamma(\ZZ\lex \RR,(1,0))$.

\begin{theorem}\label{theorem:9}
For every cardinal $\alpha$, there exists an ultrapower $[0,1]_{\rm lex}^{*}$ of the MV-algebra $[0,1]_{\rm lex}$ 
such that  any perfect MV-algebra $A$ of cardinality smaller than or equal to $\alpha$, can be embedded 
into an MV-algebra of functions from some set $X$ to $[0,1]_{\rm lex}^{*}$.
\end{theorem}
\begin{proof}
Let $A$ be a perfect MV-algebra such that $|A|=\alpha$, with $\alpha$ an infinite cardinal. Then, By Chang representation theorem, $A$ is embeddable in the product $\prod_{P \in \Spec A}A/P$ where each $A/P$ is a perfect MV-chain. 
The \el group $\mathcal{D}(A/P)$ is a totally ordered group hence it is embeddable in 
a divisible ordered group that in turns is elementary equivalent to $\RR$. 

Let $F$ be a $\alpha$-regular ultrafilter over $\alpha$; then, by \autoref{thm:chang}, $\prod_{F}\RR$ is $\alpha^{+}$-universal, hence for every $P\in \Spec A$, $\mathcal{D(A/P)} \hookrightarrow \prod_{F}\RR$, and so
$A/P \hookrightarrow \prod_F [0,1]_{\rm lex}$ (see \autoref{th:MVultrapower}). Combining the embeddings we have that
\[A \hookrightarrow \left(\prod_{F}[0,1]_{\rm lex}\right)^{\Spec A}.\]
\end{proof}

Such a result can be further improved by showing that, for any infinite cardinal $\alpha$, there exists an ultrapower of $[0,1]_{\rm lex}$ (more precisely, an iterated ultrapower), which is definable in $\alpha$  and such that any perfect MV-algebra of cardinality at most $\alpha$ embeds in the algebra of functions $\left(\prod[0,1]_{\rm lex}\right)^{X}$, for some set $X$.

A detailed proof of the result claimed above would require a full explanation on how iterated ultrapowers are constructed.  Roughly, an iterated ultrapower is obtained by iterating the ultrapower construction using a set of linearly ordered  ultrafilters $S$. The structure resulting from such a construction has the remarkable property of containing an isomorphic copy of every ultrapower defined on an ultrafilter in $S$. The theorem below formalizes the result claimed above and gives a sketch of proof; we suggest the interested reader to consult \cite{ChK,DNLS}.

\begin{theorem}\label{theorem:10}
For any infinite cardinal $\alpha$, there exists an ultrapower of $[0,1]_{\rm lex}$ definable in $\alpha$, denoted by $\prod_{\alpha}[0,1]_{\rm lex}$,
such that any perfect MV-algebra $A$, with $|A|\leq \alpha$, can be embedded into an MV-algebra of functions from some set $X$ to $\prod_{\alpha}[0,1]_{\rm lex}$.
\end{theorem}
\begin{proof}
The idea is to construct an iterated ultrapower of $[0,1]_{\rm lex}$ using the set of \emph{all} ultrafilters on $\alpha$, obviously such a structure will have the desired properties. The only obstruction is that there is no definable linear order on the set of ultrafilter on $\alpha$, if we do not admit repetitions. Since repetitions do not affect the contruction of the iterated ultrapower, we linearly order the ultrafilters as follows.

Let $Y$ be the set of all maps $y:|{\P}(\alpha)|\rightarrow {\P}(\alpha)$ such that the image
of $y$ is an ultrafilter on $\alpha$. Note that every ultrafilter on $\alpha$ appears as image of some (actually infinitely many) elements of $Y$.

The set $Y$ is totally ordered by setting $y<y'$ if there is an ordinal $\xi<| {\P}(\alpha)|$ such that
$y|_\xi=y'|_\xi$ (that is, $y$ and $y'$ coincide on all the ordinals less than $\xi$) and $y(\xi)<y'(\xi)$ in the lexicographic order of $P(\alpha)$.

For every $y\in Y$, let $D_y$ be the ultrafilter on $\alpha$ associated to $y$, i.e$.$ the image of $y$.
Then $D_\alpha$ is the resulting indexed family of ultrafilters on $\alpha$.
\end{proof}

\section{Local MV-algebras and $\ell u$-groups}

Local MV-algebras are MV-algebras with only one maximal ideal that, hence, contains all
infinitesimal elements.
This class of algebras includes MV-chains and perfect MV-algebras and it has an important
role in the representation theory of MV-algebras.
Indeed all MV-algebras can be represented as sections over sheaves whose stalks are local MV-algebras: see for example \cite{filgeo}, where the base topological space is the space of maximal ideals, and \cite{ferlet}, where the base space is $\Spec$. We recall some basic properties of local MV-algebras.

\begin{definition}\label{order}
For any $x\in A$, the order of $x$, in symbols $\ord(x)$, is the smallest
natural number $n$ such that $nx=1$.
If no such $n$ exists, then $\ord(x)=\infty .$
\end{definition}

  Let $A$ be an MV-algebra. We recall that an element $x$ of $A$ is
called \emph{finite} iff $\ord(x)<\infty $ and $\ord(\lnot x)<\infty $. Of
course $x$ is finite iff $x\wedge \lnot x$ is finite.

\begin{definition}\label{localmv}
An MV-algebra $A$ is called \emph{local} if it has only one maximal ideal, coinciding with
$\Rad(A)=\{a \in A \mid \ord(a)=\infty\}$.
\end{definition}

If $A$ is a local MV-algebra then $A/\Rad A$ is a simple MV-algebra, since it does not have non-trivial ideals.

Recall that, for any $\ell$-group $G$, a subgroup $J$ of $G$ is called an \emph{\el ideal} if it satisfies the following additional condition:
$$\textrm{if $x \in J$ and $\left|y\right| \leq \left|x\right|$ then $y \in J$.}$$
It is well-known that \el ideals are in one-one correspondence with \el group congruences and that, if $G$ is an $\ell u$-group,
then $G/J$ is an $\ell u$-group with strong unit $u/J$.
\begin{definition}\label{locallg}
An \el group $G$ is called \emph{local} if it has only one maximal \el ideal.
\end{definition}

Perfect MV-algebras and MV-chains are local MV-algebras, while the converse is false. In fact,
if $G$ is an $\ell$-group, then
\begin{equation*}
\Gamma (\ZZ\lex G,(2,0))
\end{equation*}
is a local MV-algebra, but if $G$ is not totally ordered, it is neither an
MV-chain nor a perfect MV-algebra.

\begin{proposition}
An MV-algebra $A$ is local if and only if for every $x \in A$, either
$\ord(x)<\infty$ or $\ord(\neg x)<\infty$.
\end{proposition}

We describe an example of local MV-algebra that will result as a sort of prototypical one. Let $X$ be an arbitrary non empty set, $U$ an MV-algebra, and $\mathbf{K}(U^{X})$ the subset of the MV-algebra $U^{X}$ defined as follows:
\begin{equation*}
\mathbf{K}(U^{X})=\{f\in U^{X}\mid f(X)\subseteq a/\Rad(U)\quad
\text{for some} \ a\in U\}.
\end{equation*}
   $\mathbf{K}(U^{X})$
will be called the \emph{the full MV-algebra of quasi-constant
functions} from $X$ to $U$ and any element $f$ of
$\mathbf{K}(U^{X})$ will be called a \emph{quasi-constant}
function from $X$ to $U$.
With the above notations we have:
\begin{proposition}{\rm \cite{dinespger}}
\label{prop:esloc} $\mathbf{K}(U^{X})$ is a local MV-algebra.
\end{proposition}
\begin{proof}
Let us show that $\mathbf{K}(U^{X})$ is a subalgebra of $U^{X}$. The zero-constant function $f_0$ belongs to $\mathbf{K}(U^{X})$ because $f_0(X) = \{0\} \subseteq \Rad(U) = 0/\Rad(U)$. If $f$ satisfies $f(X) \subseteq a/\Rad(U)$, for some $a\in U$, then we have $\lnot f(X) \subseteq \lnot a/\Rad(U)$. Finally, let $f, g \in \mathbf{K}(U^{X})$ be such that $f(X) \subseteq a/\Rad(U)$ and $g(X) \subseteq b/\Rad(U)$, with $a, b\in U$. Then $(f \oplus g)(X) \subseteq (a \oplus b)/\Rad(U)$. Hence $\mathbf{K}(U^{X})$ is a subalgebra of $U^{X}$.

In order to show that $\mathbf{K}(U^{X})$ is local, let us consider $f \in \mathbf{K}(U^{X})$. If $f(X) \subseteq 0/\Rad(U)$, then $\lnot f(X) \subseteq 1/\Rad(U)$ and $\ord(\lnot f) < \infty$. If $f(X) \subseteq 1/\Rad(U)$, we have $\ord(f) < \infty$. Assume now that $f(X) \subseteq a/\Rad(U) \notin \{0/\Rad(U),1/\Rad(U)\}$. Then, for every $x \in X$, $f(x) \equiv_{\Rad(U)} a$ and $a \notin \Rad(U)$. Since $A/\Rad(U)$ is an MV-chain, then $\ord(f) < \infty$ and $\ord(\lnot f) < \infty$.
\end{proof}

\begin{lemma}{\rm \cite{dinlet}}\label{lem:simpleform}
Let $x,y\in [0,1]$, and $x<y$. Then there is a  term $\varphi $ such that $\varphi (x)=0$ and $\varphi (y)=1$.
\end{lemma}

The following result is a representation based on the spectral properties of local MV-algebras.

\begin{proposition}\label{prop:spec}
Let $A$ be a local MV-algebra. Then for every $x \in A$ and for every $P, Q \in \Spec A$, we have
\begin{equation*}
\frac{(x/P)}{\Rad(A/P)}=\frac{(x/Q)}{\Rad(A/Q)}.
\end{equation*}
\end{proposition}
\begin{proof}
We observe that for every $P \in \Spec A$, $(A/P)/\Rad(A/P)$ is simple and then, up to an isomorphism, is a subalgebra of $[0,1]$.
Assume, by contradiction, that there are $r,s\in [0,1]$, with $r<s$, such that $r=(x/P)/\Rad(A/P)$ and $s=(x/P)/\Rad(A/P)$.
By Lemma \ref{lem:simpleform}, there is a simple term $\varphi$ such that $\varphi(r)=0$ and $\varphi (s)=1$.
Then $\varphi((x/P)/\Rad(A/P))=(\varphi(x)/P)/\Rad(A/P)=0$ hence
$\varphi(x)/P \in \Rad(A/P)$ while $\varphi((x/Q)/\Rad(A/Q))=(\varphi(x)/Q)/\Rad(A/Q)=1$ hence $\neg \varphi(x)/Q \in \Rad(A/Q)$.
Thus $\ord(\varphi (x))=\infty $ and $\ord(\lnot \varphi (x))=\infty $, in contrast with the assumption that $A$ is local.
\end{proof}

\begin{theorem}{\rm \cite{dinespger}}
\label{th:radconst} Every local MV-algebra can be embedded into an MV-algebra of quasi-constant functions on an ultrapower of $[0,1]$.
\end{theorem}
\begin{proof}
Let $A$ be a local MV-algebra.
By \autoref{theo:DNLT}, any element $x$ of $A$ is a function from $\Spec A$ into $[0,1]^*$,
hence $x/P \in [0,1]^*$ for every  $P \in \Spec A$. Moreover, by Proposition \ref{prop:spec},
for any $x \in A$ there exists $r_x \in [0,1]$ such that $(x/P)/\Rad(A/P) = r_x$.
Now, since every $A/P$ is embeddable in $[0,1]^*$, $\Rad(A/P)$ is embeddable in $\Rad([0,1]^*)$.
Therefore, for every $P \in \Spec A$, we have $x/P \subseteq r_x/\Rad[0,1]^*$, whence $A$ is an algebra of quasi-constant functions.
\end{proof}

Next, we shall state and prove an analogous representation for local $\ell u$-groups.
\begin{theorem}{\rm \cite{cdm}}\label{ideals}
Let $(G, u)$ be an $\ell u$-group and $A = \Gamma(G,u)$. Then the posets $\mathcal I(G)$ of all $\ell$-ideals of $G$ and $\mathcal I(A)$ of all ideals of $A$, ordered by set inclusion, are isomorphic under the following maps:
\begin{enumerate}
\item[$\phi:$] $J \in \mathcal I(A) \longmapsto \{x \in G \mid \left|x\right| \wedge u \in J\} \in \mathcal I(G)$,
\item[$\psi:$] $H \in \mathcal I(G) \longmapsto H \cap [0,u] \in \mathcal I(A)$.
\end{enumerate}
\end{theorem}

\begin{corollary}\label{local}
An $\ell u$-group $(G, u)$ is local if and only if $\Gamma(G,u)$ is a local MV-algebra.
\end{corollary}

\begin{theorem}\label{groupquasiconstant}
Every local $\ell u$-group is embeddable in an \el group
of quasi-constant functions over an ultrapower of $\RR$.
\end{theorem}
\begin{proof}
Let $\RR$ be the totally ordered group of real numbers, $I$ an ordered set and $U$ an ultrafilter of $I$. Since the property of being a totally ordered group is a first order property, then the ultrapower $R^I/U$ is again a totally ordered group.

Let $X$ be any set and consider the set $K(X,\RR^I,U)=$
$$
\left\{f: X \to \RR^I \mid \textrm{ there exists } r \in \RR^I \textrm{ such that } f(x) \in [r]_U \ \forall x \in X\right\}.
$$
$K(X,\RR^I,U)$ is an \el group with pointwise defined operations and order. Denoting by $\underline{0}$ the identity of $\RR^I$, the element $\mathbf{0}: x \in X \lmapsto \underline{0} \in \RR^I$, is the identity of $K(X,\RR^I,U)$. Analogously, $\mathbf{1}: x \in X \lmapsto \underline{1} \in \RR^I$
is a strong unit of $K(X,\RR^I,U)$.\footnote{Note that from the ultrapower construction we are not guaranteed that there exists a strong unit.}

Then $\Gamma(K(X,\RR^I,U),{\bf 1})=$
$$
\left\{f: X \to [0,1]^I \mid \textrm{there exists } r \in \RR^I \textrm{ such that }
f(x) \in [r]_U \ \forall x \in X\right\},
$$
is precisely the MV-algebra $\mathbf{K}(([0,1]^*)^X)$.

Now, if $(G,u)$ is a local $\ell u$-group, by Corollary \ref{local}, $A = \Gamma(G,u)$ is a local MV-algebra. Hence, by \autoref{th:radconst}, $A$ embeds in some $\mathbf{K}(([0,1]^*)^X)$. So, since $\Gamma$ and $\Xi$ are equivalences, $(G,u)$ embeds in $\Xi(\mathbf{K}(([0,1]^*)^X)) = K(X,\RR^I,U)$.
\end{proof}

The machinery used at the end of the previous section for perfect MV-algebras can also be used here to find uniform and canonical representation.  We state the theorem without proof as it is an easy adaptation of the proofs of \autoref{theorem:9} and \autoref{theorem:10}.


\begin{theorem}
\label{uniformlocal}
For any infinite cardinal $\alpha$, there exists an ultrapower of $\RR$, definable in $\alpha$, which we call $\prod\RR_{\alpha}$, 
such that any $\ell u$-group of cardinality at most $\alpha$ is embeddable in the 
$\ell u$-group of quasi-constant functions over $\prod\RR_{\alpha}$.
\end{theorem}

\end{document}